\documentclass[a4paper,twoside,10pt,reqno]{article}
\usepackage[english]{babel}
\usepackage[dvips]{graphicx}
\usepackage[T1]{fontenc}
\usepackage[latin1]{inputenc}
\usepackage{amssymb,amsthm,dsfont,mathrsfs}
\usepackage{amsfonts,amsmath,euscript}
\usepackage{color,cite}

\newtheorem{teo}{Theorem}[section]
\newtheorem{pro}[teo]{Proposition}

\newtheorem{ex}{Example}[section]

\newcounter{example}[section]

           \addto{\captionsenglish}{}

\newcommand{\hs}{\hspace{3pt}}
\newcommand{\dst}{\displaystyle}
\newcommand{\dem}{{\bf Dem. }}
\newcommand{\fdem}{$\square$}

\newcommand{\titulo}[1]{\mbox{} \\ \noindent \textit{\textbf{\Large #1}}\\}

\renewcommand{\abstract}[1]{{\small \noindent \textbf{Abstract:} #1\\}}

\newcommand{\pchave}[1]{{\small \noindent \textbf{Keywords:} #1\\}}

\begin{document}

\begin{center}
\titulo{Generating multivariate extreme value distributions}
\end{center}

\vspace{0.5cm}

\textbf{Helena Ferreira} Department of Mathematics, University of Beira
Interior, Covilhã, Portugal\\

\abstract{We define in a probabilistic way a parametric family of multivariate extreme value distributions. We derive its copula, which is a mixture of several complete dependent copulas and total independent copulas, and the bivariate tail dependence and extremal coefficients. Based on the obtained results for these coefficients, we propose a method to built multivariate extreme value distributions with prescribed tail/extremal coefficients. We illustrate the results with examples of simulation of these distributions.}

\pchave{multivariate extreme value theory, tail dependence, extremal
coefficients, simulation}

\section{Introduction}

The construction of multivariate distributions has its motivation in Probability Theory, Biostatistics, Economics, and there is no need to demonstrate its importance nowadays. In particular, multivariate extreme value (MEV) distributions provide models for joint extreme events which require special care from practitioners since rare events can have serious natural and economic impact.\\
The MEV distributions composes an important class of positive dependent distributions with dependence in the extreme upper tails, with this property known as tail dependence (\cite{resnick, joe}). A misleading evaluation of the dependence in the tails of a multivariate distribution may lead to an underestimation of risks from crises events (\cite{emb+02}).\\
The MEV distributions can be constructed through methods such as mixtures, stochastic representations and limits (\cite{joe}). They also can be obtained via copulas and suitable techniques of extra-parametrisation (\cite{nelsen_cop,fisher,liebscher,ferreira+pereira}).\\ 
The dependence structure of a MEV distribution is completely characterised by its dependence function (\cite{resnick}). However this function cannot be easily inferred from data and simple dependence measures are welcome.\\
The most popular of the dependence measures for bivariate distributions are the upper tail dependence coefficient $\lambda$ and the extremal coefficient $\epsilon$, introduced far back in the sixties (\cite{sib,tiago}), which determines each other in the extreme value distributions by the relation $\epsilon=2-\lambda$.\\
Our paper is a contribution for the main issue of the development of higher dimensional models capturing dependencies. This issue covers quite a large spectrum of techniques and applications where copulas become ladies of quality since they are able to yield any kind of dependence structure independently of marginal distributions. An application of copula functions is the simulation of multivariate distributions with dependent observations. By departing of the right copula or a mixture of copulas any dependence structure may be reached. For a review of different estimators, some model selection tests for copula functions and methods of its simulation see, for instance,  \cite{ferm, gennest,patton}. However, in general copulas cannot be characterised in a simple way that provides a simple algorithm to simulate data.\\
How to get out of the way of selection of a candidate copula, its estimation and simulation of dependent random observations? One proposal is to estimate non-parametricaly the bivariate extremal dependences (\cite{schmidt+stadt, ferreiras}) and to construct in a probabilistic way a multivariate extreme value distribution with these prescribed extremal dependences, which can be used to generate random observations.\\
Here we propose a method to built a MEV distribution with prescribed tail/extremal dependence coefficients.\\
To the best of our knowledge, the problem was treated by \cite{brommundt} for elliptical distributions and by \cite{falk} for the case where each bivariate tail dependence coefficient $\lambda_{ij}$ between the margins $X_i$ and $X_j$ has the representation $\lambda_{ij}=2-(1+a_i-a_j)^{1/\alpha}$, where $0\leq a_1\leq ...\leq a_d\leq 1\leq \alpha$ and $d$ is the dimension of the random vector.\\
Our solution is of stochastic representation kind and can be used in practical problems where one needs to buil a stochastic model in a situation where the extremal dependence is known or estimated and the knowledge or estimates of marginal distributions are available.\\
 A motivation for such approach can be the difficulty to choose or find an appropriate copula for the problem in hand or the fact that the choice of copulas does not inform explicitly about the strength of the dependence between the variables involved.
Another advantage of the method we propose is that we can take tail dependence/extremal coefficients not constant over  time and build a more realistic time-varying parametric family of models.\\
The result shows that the MEV distributions are rich enough to encompass a large variability for bivariate dependence one would like to be able to handle.

\section{A multivariate extreme value model}\label{sed}

Let $Y_i$, $i=1,...,d\geq 2$, and $Z_j$, $j=1,...,D$, be independent and unit Frech\'{e}t variables and $\alpha_j^{(i)}$, $i=1,...,d, j=1,...,D,$  nonnegative constants. For $C$ constant satisfying $C\geq \displaystyle\bigvee_{i=1}^d\displaystyle\sum_{j=1}^D\alpha_j^{(i)}$,  consider the $d$-dimensional random vector ${\bf X}=(X_1,...,X_d)$ defined as 
\begin{eqnarray}\label{1}
X_i=\bigvee_{j=1}^{D}\left(\alpha_j^{(i)} Z_j\right)\bigvee \left ( C-\sum_{j=1}^{D}\alpha_j^{(i)}\right )Y_i,\,\,\,i=1,...,d.
\end{eqnarray}

We will present in the next result the distribution of the random vector ${\bf X}$ and its copula  $C_{{\bf X}}(u_1,...,u_d)=F_{{\bf X}}(F^{-1}(u_1),...,F^{-1}(u_d))$, where $F^{-1}$ denotes the generalised inverse of the distribution function of $X_i$, which is a mixture of several complete dependent copulas and total independent copulas.

\begin{pro}\label{p.1} The mixture model in (\ref{1}) has multivariate extreme value distribution defined as 
\begin{eqnarray}\label{2}
F_{{\bf X}}(x_1,...,x_d)=\displaystyle\prod_{j=1}^{D}\left(\bigwedge_{i=1}^{d}e^{-x_i^{-1}\alpha_j^{(i)}}\right)
\displaystyle\prod_{i=1}^{d}e^{-x_i^{-1}\left(C-\displaystyle\sum_{j=1}^{D}\alpha_j^{(i)}\right)},
\end{eqnarray}
with its corresponding copula 
\begin{eqnarray}\label{3}
C_{{\bf X}}(u_1,...,u_d)=\displaystyle\prod_{j=1}^{D}\left(\bigwedge_{i=1}^{d}u_i^{\alpha_j^{(i)}/C}\right)
\displaystyle\prod_{i=1}^{d}u_i^{\left(1-\frac{1}{C}\displaystyle\sum_{j=1}^{D}\alpha_j^{(i)}\right)}.
\end{eqnarray}
\end{pro}
\dem By assuming $1/0=+\infty$, we can write from the definition of the model in (\ref{1}), for each $x_i\geq 0, i=1,...,d$,
 \begin{eqnarray}\label{}
F_{{\bf X}}(x_1,...,x_d)=P\left(\displaystyle\bigcap_{j=1}^{D}\left\{Z_j\leq\displaystyle\bigwedge_{i=1}^{d}\frac{x_i}{\alpha_j^{(i)}}\right\},
\displaystyle\bigcap_{i=1}^{d}\left\{Y_i\leq\frac{x_i}{C-\displaystyle\sum_{j=1}^{D}\alpha_j^{(i)}}\right\}\right).
\end{eqnarray}
The assumptions on the distributions of the variables $Z_j$ and $Y_i$ leads to the expression in (\ref{2}) and, by changing to the variable $u_i=e^{-x_i^{-1}C}$ we obtain the copula function which satisfies the max-statbility condition $C_{{\bf X}}^t(u_1,...,u_d)=C_{{\bf X}}(u_1^t,...,u_d^t)$, for $t>0$. 
 \fdem\\
 
 The tail dependence coefficient between $X_i$ and $X_j$, $i,j\in \{1,...,d\}$ measures the probability of occurring extreme values
 for one random variable given that another assumes an extreme value too and is defined as 
\begin{eqnarray}\label{tdc}
\dst\lambda_{ij}=\lim_{u\uparrow 1}P(F_{X_i}(X_i)>u|F_{X_j}(X_j)>u),
\end{eqnarray}
where $F_{X_i}$ denotes the distribution function of $X_i$.
In our model, since $(X_i,X_j)$ has bivariate extreme value distribution, $\lambda_{ij}$ 
can be computed from the extremal coefficient $\epsilon_{ij}$ of  $(X_i,X_j)$. This is defined by 
\begin{eqnarray}\label{ec}
\dst\epsilon_{ij}=-\log F_{(X_i,X_j)}(x,x), x>0,
\end{eqnarray}
where $F_{(X_i,X_j)}$ denotes the distribution function of $(X_i,X_j)$. The value of $\epsilon_{ij}$ doesn't depend on $x$ and satisfies $\epsilon_{ij}=2-\lambda_{ij}$.

These dependence measures have extensions for $d$-dimensional vectors with $d>2$ ( \cite{schmidt+stadt, li3, ferreiras, schlather+tawn}).\\
In the next result we compute $\lambda_{ij}$ for the model in (\ref{1}) and show how such a model for prescribed $\lambda_{ij}$ can be constructed. In this way we provide a method for building new multivariate distributions and modelling dependence structures.

\begin{pro}\label{p.2} (a) The random vector ${\bf X}$ defined  in (\ref{1}) has tail dependence coefficients 
\begin{eqnarray}\label{4}
\dst\lambda_{sk}=\frac{1}{C}\displaystyle\sum_{j=1}^{D}\left(\alpha_j^{(s)}\wedge\alpha_j^{(k)}\right),\,\,\,s,k\in\{1,...,d\},\,\,\,s\neq k.
\end{eqnarray}
(b) Given the set of tail dependence coefficients $\{\lambda_{sk},\,\,s,k\in\{1,...,d\},\,\,s\neq k\}$ of a $d$-dimensional random vector and a constant $C\geq d-1$, there exists a random vector ${\bf X}$ defined as  in (\ref{1}) with proportional tail dependence coefficients $\{\frac{1}{C}\lambda_{sk},\,\,s,k\in\{1,...,d\},\,\,s\neq k\}$.
\end{pro}
\dem First we derive the diagonal of the bivariate copula $C_{sk}(u,v)$ of $(X_s,X_k)$, $1\leq s<k\leq d$. We have, by taking $u_i=1$ for $i\not\in\{s,k\}$ in (\ref{3}), 
\begin{eqnarray}
C_{sk}(u,u)=u^{\frac{1}{C}\displaystyle\sum_{j=1}^{D}\left(\alpha_j^{(s)}\vee\alpha_j^{(k)}\right)}u^{1-\frac{1}{C}\displaystyle\sum_{j=1}^{D}\alpha_j^{(s)}}u^{1-\frac{1}{C}\displaystyle\sum_{j=1}^{D}\alpha_j^{(k)}}=
u^{2-\frac{1}{C}\displaystyle\sum_{j=1}^{D}\left(\alpha_j^{(s)}\wedge\alpha_j^{(k)}\right)}.
\end{eqnarray}
Therefore 
\begin{eqnarray}
\epsilon_{ks}=\epsilon_{sk}=-\log F_{(X_s,X_k)}(C,C)=-\log C_{sk}(e^{-1},e^{-1})=2-\frac{1}{C}\displaystyle\sum_{j=1}^{D}\left(\alpha_j^{(s)}\wedge\alpha_j^{(k)}\right),
\end{eqnarray}
which leads to the above presented value for $\lambda_{sk}$.

Now, for (b), let $\Lambda=[\lambda_{sk}]_{s,k=1,...,d}$ be the symmetric tail dependence coefficients matrix of a random vector, where $\lambda_{ss}=1,\,\,s=1,...,d$. Consider a model as in (\ref{1}) with $\alpha_j^{(i)}$, $i=1,...,d$, $j=1,...,D=\frac{d(d-1)}{2}$ defined as following.
Assume that when $j\in \emptyset$ we don't define the corresponding coefficient,  $n(0)=0$, $n(i)=n(i-1)+d-i,\,\, i\in\{1,...,d-1\}$ and $n(d)=\frac{d(d-1)}{2}$.  Let us consider
\begin{eqnarray}\label{5}
\alpha_j^{(i)}=\left\{
\begin{array}{rcl}
m_1\delta_{j-i,n(0)-1} & \mbox{if} & j\in \{n(0)+1,...,n(1)\}\\
m_2\delta_{j-i,n(1)-2} & \mbox{if} & j\in \{n(1)+1,...,n(2)\}\\
...\\
m_{i-1}\delta_{j-i,n(i-1)-(i-1)} & \mbox{if} & j\in \{n(i-2)+1,...,n(i-1)\},
\end{array}
\right.
\end{eqnarray}
\begin{eqnarray}\label{6}
\alpha_j^{(i)}=\lambda_{i,j+i-n(i-1)},\,\,\, j\in \{n(i-1)+1,...,n(i)\},
\end{eqnarray}
\begin{eqnarray}\label{7}
\alpha_j^{(i)}=0,\,\,\, j\in \left\{n(i)+1,...,\frac{d(d-1)}{2}\right\},
\end{eqnarray}

with $\delta_{i,j}=1$ if $i=j$,  $\delta_{i,j}=0$ if $i\neq j$ and $m_i=\displaystyle\bigvee_{j=i+1}^d\lambda_{ij}$.

We remark that $n(d-1)=\frac{d(d-1)}{2}$ and therefore the set in (\ref{7}) will be empty for $ i=d-1$. Also, since $n(0)=0$, the sets in (\ref{5}) will be empty for $i=1$.\\
For sake of simplicity we can present the coefficients $\alpha_j^{(i)}$ in a matrix $A=[\alpha_j^{(i)}]_{i=1,...,d, j=1,...,D}$ as

$$A=\left [
\begin{array}{ccccccccccccccccc}

\lambda_{12} & \lambda_{13} & ....&... & \lambda_{1d} & 0 & &...& &0 & 0&...&...  &...& &...  &0\\
m_1& 0&...&0 &0&\lambda_{23} & ...&...& & \lambda_{2d}   &0&...&...&...&... & &0\\
0& m_1& ...& 0&0&  m_2& 0&...& 0& 0 &\lambda_{34} & ...&\lambda_{3d} &  0 &...& &0\\
...&...&...&...&....&0&m_2& &  & 0&m_3&...&0& ...&...& &...\\
...&...&...&...&....& & & &   & 0&...&...&...& ...&...& &...\\
0&...&     0&m_1&0&0&...&   & &0&   &...&...&... & & ...&                     \lambda_{d-1,d}\\
0&...&     0&0&m_1&0&  & &...&m_2 &0&  ...&m_3&...  & ...& &   m_{d-1}\\
\end{array}
\right].$$

\vspace{0,5cm}
By assuming $\sum_{j=p}^q a_j=0$ if $q<p$, we can verify that  this choice of the coefficients for (\ref{1}) leads to  
\begin{eqnarray}\label{7}
\displaystyle\sum_{j=1}^D\alpha_j^{(i)}=\displaystyle\sum_{j=1}^{i-1}m_j+\displaystyle\sum_{j=i+1}^d\lambda_{ij}\leq i-1+d-i=d-1\leq C,\,\,\,i=1,...,d,
\end{eqnarray}

and, for each pair $1\leq s <k\leq d$,
 
$$\frac{1}{C}\displaystyle\sum_{j=1}^{D}\left(\alpha_j^{(s)}\wedge\alpha_j^{(k)}\right)=
\frac{1}{C}\displaystyle\sum_{j=1}^{n(s-1)}\left(\alpha_j^{(s)}\wedge\alpha_j^{(k)}\right)+
\frac{1}{C}\displaystyle\sum_{j=n(s-1)+1}^{n(s)}\left(\alpha_j^{(s)}\wedge\alpha_j^{(k)}\right)+
\frac{1}{C}\displaystyle\sum_{j=n(s)+1}^{\frac{d(d-1)}{2}}\left(\alpha_j^{(s)}\wedge\alpha_j^{(k)}\right)=$$
$$0+\frac{1}{C}(\lambda_{sk}\wedge m_s)+0=\frac{1}{C}\lambda_{sk}. $$ 

\hspace{15cm}\fdem\\

 In fact, as we can seen in the proof, the model can be constructed for any 
 $$C\geq \displaystyle\bigvee_{i=1}^d \left(\displaystyle\sum_{j=1}^{i-1}m_j+\displaystyle\sum_{j=i+1}^d\lambda_{ij}\right)$$
 and therefore we have the next two particular cases.\\

 \begin{pro}\label{p.3} (a) Given the set of tail dependence coefficients $\{\lambda_{sk},\,\,s,k\in\{1,...,d\},\,\,s\neq k\}$ of a $d$-dimensional random vector such that 
 $\displaystyle\sum_{j=1}^{i-1}\left(\displaystyle\bigvee_{k=j+1}^d\lambda_{jk}\right)+\displaystyle\sum_{j=i+1}^d\lambda_{ij}\leq 1$, for each $i=1,...,d$,  there exists a random vector ${\bf X}$ defined as  in (\ref{1}) with the same tail dependence coefficients.\\
 
\noindent
(b) Given the set of tail dependence coefficients $\{\lambda_{sk},\,\,s,k\in\{1,...,d\},\,\,s\neq k\}$ of a  $d$-dimensional random vector such that $\lambda_{sk}\leq \frac{1}{d-1}$ for each $1\leq s <k\leq d$, there exists a random vector ${\bf X}$ defined as  in (\ref{1}) with the same tail dependence coefficients.
\end{pro}

 \section{Generating MEV distributions}

\begin{ex} Consider the model in (\ref{1}) with $\alpha_1^{(1)}=\frac{1}{2}$,  $\alpha_2^{(1)}=2$,   
$\alpha_1^{(2)}=\frac{1}{4}$, $\alpha_2^{(2)}=2$, 
$\alpha_1^{(3)}=1$, $\alpha_2^{(3)}=\frac{1}{2}$ and $C=\frac{5}{2}$. By applying the Proposition 1.2, we find the following tail dependence matrix
$$\Lambda=\left [
\begin{array}{ccc}

1 & 0.9&0.4\\
0.9& 1&0.3\\
0.4& 0.3& 1\\
\end{array}
\right]$$

\noindent
for the vector ${\bf X}=(X_1,X_2,X_3)$.

\end{ex}
\begin{ex} Let now 
$$\Lambda=\left [
\begin{array}{cccc}

1 & 0.2&0.5&0.3\\
0.2& 1&0.6&0.1\\
0.5& 0.6& 1&0.9\\
0.3&0.1&  0.9&1\\
\end{array}
\right]$$

\noindent
be the tail dependence coefficients matrix of a $4$-dimensional random vector.\\
We chose 
$$C=(0.2+0.5+0.3)\vee\left((0.2\vee 0.5\vee 0.3)+0.6+0.1\right)\vee \left((0.2\vee 0.5\vee 0.3)+(0.6\vee 0.1)+0.9\right)=2$$
and, as in the proof of part (b) of the Proposition 1.2, we present the coefficients $\alpha_j^{(i)}$, $j=1,...,6$, for $X_i$ along line $i$ of  the matrix $A$, $i=1,...,4$, as follows
$$A=\left [
\begin{array}{ccccccc}

0.2 & 0.5&0.3&0&0&0\\
0.5& 0&0&0.6&0.1&0\\
0& 0.5& 0&0.6&0&0.9\\
0&0&0.5&0&0.6&0.9\\
\end{array}
\right].$$

\noindent
For this choice of $\alpha_j^{(i)}$, $i=1,...,4,\,\,j=1,...,6$ and $C$, by applying  the Proposition 1.1, we obtain a random vector ${\bf X}=(X_1,X_2,X_3,X_4)$ with tail dependence coefficients matrix

$$\Lambda^{*}=\left [
\begin{array}{cccc}

1 & 0.1&0.25&0.15\\
0.1& 1&0.3&0.05\\
0.25& 0.3& 1&0.45\\
0.15&0.05&  0.45&1\\
\end{array}
\right].$$
\end{ex}

\begin{ex} Let now 
$$\Lambda=\left [
\begin{array}{ccc}
1 & 0.2&0.1\\
0.2& 1&0.8\\
0.1& 0.8& 1\\
\end{array}
\right]$$

\noindent
be the tail dependence coefficients matrix of a $3$-dimensional random vector.\\
If we chose $C=(0.2+0.1)\vee\left((0.2\vee 0.1)+0.8\right)=1$, the matrix of coefficients $\alpha_j^{(i)}$, $i=1,...,3,\,\,j=1,...,3$,
$$A=\left [
\begin{array}{ccc}
0.2 & 0.1&0\\
0.2& 0&0.8\\
0& 0.2& 0.8\\
\end{array}
\right]$$
then, by applying Proposition 2.3,  the model in (\ref{1}) still has the given tail dependence coefficients matrix $\Lambda$.\\

\noindent
We simulated $1000$ observations  of ${\bf X}$ leading to the following bivariate representations.
\begin{figure}
\begin{center}
\includegraphics[clip,width=12cm]{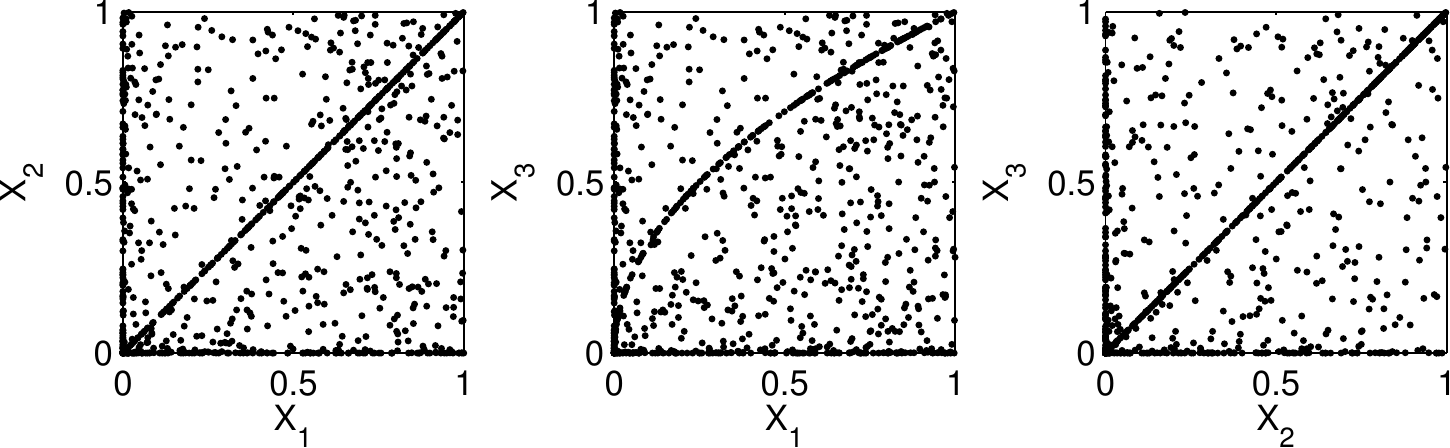}
\caption{Simulated values for:  $(X_1,X_2)$ (left), $(X_1,X_3)$ (center) and $(X_2,X_3)$ (right).}
\end{center}
\end{figure}

\end{ex}

\end{document}